\documentclass{article}
\usepackage{arxiv}
\usepackage{graphicx}
\usepackage{float}
\usepackage{amsmath,amssymb}
\usepackage{amsthm}
\usepackage{algorithm,algorithmic}
\usepackage[labelformat=simple]{subfig}

\usepackage{multirow}
\usepackage{latexsym,bm}
\usepackage{hyperref}
\usepackage{doi}     

\newcommand {\lr}[1] {\left\langle #1 \right\rangle_{\Omega}}
\newcommand {\lrh}[1] {\left\langle #1 \right\rangle_{\mathcal{T}_h}}
\graphicspath{{./fig/}}

\usepackage{hyperref}
\makeindex

\title{Reduced order modeling for spatio-temporal pattern approximation in diffusive Lotka-Volterra equations}

\stitle{ROM for spatio-temporal pattern approximation in diffusive Lotka-Volterra equations}

\author{
  B\"ulent Karas\"ozen\\
     Department of Mathematics \& Institute of Applied Mathematics\\
     Middle East Technical University\\
     Ankara-Turkey\\
     \texttt{bulent@metu.edu.tr}\\
  \And
  Murat Uzunca\\
   Department of Mathematics\\ Sinop University\\
     Sinop-Turkey\\ \
     \texttt{muzunca@sinop.edu.tr}\\
}

\begin{document}

\maketitle

\begin{abstract}
This paper presents an efficient reduced order modeling (ROM) framework for simulating spatio-temporal pattern formation in three-species diffusive Lotka-Volterra systems. To alleviate the high computational cost associated with long-time simulations of the high-dimensional full order model (FOM), we apply proper orthogonal decomposition (POD) to project the solution onto a low-dimensional subspace. Further efficiency is achieved through tensorial POD (TPOD), which preserves the quadratic nonlinear structure and enables offline-online decomposition. Numerical experiments demonstrate that both POD and TPOD accurately replicate the key features of spatial segregation patterns while substantially reducing computation time, whereas  the TPOD is faster.  Additionally, we demonstrate accurate long-time pattern prediction  using limited training data.
\end{abstract}



\section{Introduction}

A central problem in population dynamics is understanding of interaction between biological species. 
The dynamics of populations are represented in terms of birth, death, and movement of the species. 
The most popular model is the Lotka-Volterra system (LVS) with the quadratic reaction terms \cite{Adamson12,Contento15,Manna21,Mimura15,Petrovskii01}.
In  this work, we consider on a two dimensional bounded domain $\Omega \in \mathbb{R}^2$, and for a time period $[0,T]$ ($T>0$), the following three species LVS  
\begin{align}\label{lv3}
\partial_t u_i & = d_i \Delta u_i  +  r_iu_i - \sum_{j=1}^3 b_{ij}u_iu_j, & i=1,2,3,
\end{align}
where the unknown functions $u_i(\mathrm{x},t)$ are population densities,  $d_i>0$ are the diffusion rates, $r_i$ the intrinsic growth rate of the $i$-th species, $b_{ii}>0$ is the intra-specific competition rate between the same species,  $b_{ij}>0$ ($i\neq j$) is the inter-specific interaction rate between different species.
The operators $\partial_t$ and $\Delta $ denote the differentiation with respect to time variable $t$ and the Laplace operator, respectively.
No-flux boundary conditions are assumed, i.e., 
$\partial_{\nu} u_i (\mathrm{x},t)=0$, $\mathrm{x} \in \partial\Omega, \; t \in [0,T]$, with the initial conditions  
$
u_i(\mathrm{x},0) \ge 0, \; \mathrm{x} \in \Omega
$,
are prescribed.

Various scenario of spiral waves and patterns in two-dimensions are reported in \cite{Adamson12,Contento15,Mimura15,Cangiani16,Ei99}. 
In order to construct the different patterns accurately, large number of numerical simulations have to be performed for long term computation.
This paper aims to apply intrusive model order reduction (MOR) to reduce the computing  time of solving numerically LVS \eqref{lv3} while maintaining the essential characteristics of the pattern-forming dynamics.  
We apply the proper orthogonal decomposition (POD) \cite{Berkooz93,Sirovich87}, which is the most known intrusive MOR method \cite{Benner17}.
Recently, MOR has been applied to pattern forming systems in \cite{Alla23,Bozzini21,Karasozen17} for Schnakenberg, and FitzHugh-Nagumo equations in \cite{Karasozen21}.

Several numerical methods like finite differences and finite elements, have been applied to solve population systems. 
We discretize LVS \eqref{lv3} in space by the spectral element method (SEM) which  is based on the Galerkin approach. The SEM yields highly accurate  solutions using a coarse mesh  \cite{Canuto06, Canuto07}. 
The semi-discrete system of ordinary differential equations (ODE) is discretized in time with linearly implicit  Kahan's method \cite{Kahan93,Celledoni13}, 
which has an algorithmic complexity that lies between fully-implicit
and semi-implicit time-stepping schemes. 
In this way, accurate and fast full-order solutions are obtained.
An efficient MOR technique is essential, as long-time simulations of high-dimensional ODE problems often render computations using the computationally prohibitive standard POD method.
In this paper, we utilize tensor techniques \cite{Karasozen21,Benner21,Benner18} for further acceleration in computation of the reduced solutions, which is suitable for linear-quadratic systems. 
We show that speedup factor of POD in tensor form is much higher than the standard POD.  
We demonstrate that the patterns can be predicted for large times accurately using large set of solutions as in \cite{Karasozen21,Alla24}.

The paper is organized as follows: in Section \ref{sec:fom}, the FOM for  the LVS \eqref{lv3} is presented. 
The standard POD and POD in tensor form are given in Section \ref{sec:rom}. 
The performance and prediction capabilities of the ROMs are illustrated for LVS with segregating patterns in Section \ref{sec:num}. The paper ends with some conclusions in Section \ref{sec:con}.

\section{Full order model} 
\label{sec:fom}

We first introduce full discrete solution of the LVS \eqref{lv3}, which we call FOM throughout the paper.
 
\subsection{Space discretization}

Among the spatial discretization methods, SEM  \cite{Canuto06,Canuto07} yields an exponential rate of convergence
for smooth solutions. As compared with a conventional finite differences,  SEM require
weaker smoothness.  Like the finite elements methods, SEM is based on the Galerkin approach, which uses the weak form of the given partial differential equation (PDE). 
The weak solution  of the LVS \eqref{lv3} is given as 
\begin{equation}\label{lv2}
\lr{ \partial_t u_i,v} = d_ia(u_i,v) + r_i\lr{u_i,v} - \sum_{j=1}^3 b_{ij}\lr{u_iu_j,v} , \qquad i=1,2,3,
\end{equation}
where $\lr{\cdot, \cdot}$ denotes the usual $L^2$-inner product on the domain $\Omega$, and $a(u,v) = \lr{\nabla u, \nabla v}$ is the bilinear form. In order to construct the finite dimensional formulation, let
$ \mathcal{T}_h$ be a family of conformal, regular, and quasi-uniform partitions of the domain $ \Omega \subset \mathbb{R}^2$ into $n_e$ quadrilaterals. We introduce the following finite dimensional solution space 
$$
V_h = \{ v \in C^0 (\overline{\Omega}) : v_{|T_\ell} \in \mathbb{Q}_p(T_\ell) , \; \forall T_\ell \in \mathcal{T}_h\},
$$
where $\mathbb{Q}_p(T_\ell)$ is the space of piecewise continuous polynomials of degree at most $p\geq 1$ on the given quadrilateral $T_\ell\in\mathcal{T}_h$.
Then, the finite dimensional problem for \eqref{lv2} reads as: for almost all $t$, find $u_{i,h}(\mathrm{x},t) \in V_h$ such that $\forall v_{i,h}(\mathrm{x})\in V_h$ there holds for $i=1,2,3$
\begin{equation}\label{lvweak}
\lrh{ \partial_t u_{i,h},v_{i,h}} = d_ia_h(u_{i,h},v_{i,h}) + r_i\lrh{u_{i,h},v_{i,h}} - \sum_{j=1}^3 a_{ij}\lrh{u_{i,h}u_{j,h},v_{i,h}} ,
\end{equation}
where the broken inner-product and the broken bilinear form are given respectively by
\begin{equation*}
\lrh{u,v} = \sum_{\ell=1}^{n_e}\langle u,v\rangle_{T_\ell},\qquad
a_h(u,v) = \sum_{\ell=1}^{n_e} \langle \nabla u, \nabla v\rangle_{T_\ell}.
\end{equation*}

Typically, when using SEM, the integral computations in \eqref{lvweak} are done by the composite Legendre-Gauss-Lobatto (LGL) quadrature rule defined on the reference interval $[-1,1]$. 
One starts with the $(p + 1)$ LGL quadrature nodes and quadrature weights on the reference interval, then they are mapped into the elements $T_{\ell} \in\mathcal{T}_h$ through an affine map and tensor. Let $\mathrm{x}_{\ell , q}$ and $w_{\ell , q}$ denote the quadrature nodes and weights on a generic element $T_{\ell} \in\mathcal{T}_h$, $l = 1,..., n_e$, $q = 1,..., (p+1)^2$ (there are $(p+1)$ nodes in each space direction). For a given partition $ \mathcal{T}_h$, we use the composite LGL quadrature formula
\begin{equation}\label{eq3}
	\lrh{u(\mathrm{x}),v(\mathrm{x})} = \sum_{\ell = 1}^{n_e}\langle u(\mathrm{x}),v(\mathrm{x})\rangle_{T_\ell} \approx \sum_{\ell = 1}^{n_e} \left[\sum_{q=1}^{(p+1)^2} w_{\ell ,q}u(\mathrm{x}_{\ell, q}) v(\mathrm{x}_{\ell,q}) \right].
\end{equation}

Let us set $N = N(p)$ as the total number of non-repeated LGL quadrature nodes on $\mathcal{T}_h$, and denote by $\{\mathrm{x}_1,\ldots , \mathrm{x}_N\}$ the corresponding non-repeated LGL quadrature nodes. We can rewrite the LGL quadrature formula \eqref{eq3}
 in the following form
\begin{equation}\label{lglq}
	\lrh{u(\mathrm{x}),v(\mathrm{x})} \approx \sum_{q = 1}^{N} w_{q}u(\mathrm{x}_{q}) v(\mathrm{x}_{q}) ,
\end{equation}
where the quadrature weights $w_{q}$ accounts the non-repeated entries.
We use the quadrature nodes $\{\mathrm{x}_1,\ldots , \mathrm{x}_N\}$ as the degrees of freedom together with the nodal Lagrange basis functions $\varphi_i (\mathrm{x})$, $i=1,\ldots ,N$. Then, we obtain the semi-discrete solutions in the following form
\begin{equation}\label{sol}
u_{h,i}(\mathrm{x},t) = \sum_{j=1}^N \mathbf{u}_{j,i}(t) \varphi_i(\mathrm{x}),
\end{equation}
where $\mathbf{u}_i(t):=(\mathbf{u}_{1,i}(t),\ldots,\mathbf{u}_{N,i}(t) )^T$ is the time-dependent vector of nodal coefficients for the component $u_{h,i}$, $i=1,2 ,3$, given $\mathbf{u}_{j,i}(t)\approx u_{h,i}(\mathrm{x}_j,t)$. Next, we impose the solution \eqref{sol} into the system \eqref{lvweak}, we choose $v_{i,h}=\varphi_s(\mathrm{x})$, $s=1,\ldots ,N$, and we obtain the linear-quadratic ODEs of the form
\begin{equation}\label{odes1}
M\dot{\mathbf{u}}_i = d_iK\mathbf{u}_i + r_iM\mathbf{u}_i - \sum_{j=1}^3 b_{ij} M(\mathbf{u}_i\circ\mathbf{u}_j),\qquad i=1,\ldots ,3,
\end{equation} 
where $\circ$ denotes the element-wise (Hadamard) product of vectors. Here, $M\in\mathbb{R}^{N\times N}$ and $K\in\mathbb{R}^{N\times N}$ are the SEM mass and stiffness matrices, respectively, whose entries are given by
$$
M_{st} = \lrh{\varphi_t,\varphi_s} \; , \quad K_{st} = a_h(\varphi_t,\varphi_s).
$$
The resulting mass matrix $M$ is diagonal
$$
M_{st} = \lrh{\varphi_t,\varphi_s} \approx \sum_{q = 1}^{N} w_{q}\varphi_t(\mathrm{x}_q)\varphi_s(\mathrm{x}_q)=\left\{
\begin{matrix}
w_s, & s=t\\
0, & s\neq t
\end{matrix}\right.,
$$
due to the orthogonality of the Lagrange basis functions.
By the same features, the nonlinear terms in \eqref{lvweak} 
are given as
\begin{align*}
\lrh{u_{i,h}u_{j,h},\varphi_s} &= \lrh{ \left(\sum_{m=1}^N  \mathbf{u}_{m,i}(t) \varphi_m(\mathrm{x})\right) \left(\sum_{n=1}^N  \mathbf{u}_{n,j}(t) \varphi_n(\mathrm{x})\right),\varphi_s(\mathrm{x})} \\
&\approx \sum_{q = 1}^{N} w_{q} \left(\sum_{m=1}^N  \mathbf{u}_{m,i}(t) \varphi_m(\mathrm{x}_q)\right) \left(\sum_{n=1}^N  \mathbf{u}_{n,j}(t) \varphi_n(\mathrm{x}_q)\right)\varphi_s(\mathrm{x}_q)\\
&= M_{ss}\mathbf{u}_{s,i}(t)  \mathbf{u}_{s,j}(t).
\end{align*}

Multiplying both sides by $M^{-1}$ (which is very cheap since the mass matrix is diagonal), the system of ODEs \eqref{odes1} can be written as 
\begin{align}\label{odes0}
\dot{\mathbf{u}}_i &=  A_i\mathbf{u}_i - \sum_{j=1}^3 b_{ij}  (\mathbf{u}_i\circ\mathbf{u}_j),& i=1,2,3, 
\end{align} 
where $A_i = d_iM^{-1}K + r_iI_N$ with $I_N$ is the $N$-dimensional identity matrix. 
Let us define the $3N$-dimensional vectors 
\begin{equation*}
\mathbf{u}(t) = 
\begin{pmatrix}
\mathbf{u}_1(t) \\
\mathbf{u}_2(t) \\
\mathbf{u}_3(t)
\end{pmatrix}, \qquad 
\mathbf{u}_{(j)}(t) = 
\begin{pmatrix}
\mathbf{u}_j(t) \\
\mathbf{u}_j(t) \\
\mathbf{u}_j(t)
\end{pmatrix}, \qquad j=1,2,3,
\end{equation*} 
and the $(3N\times 3N)$-dimensional matrices
\begin{equation*}
A = 
\begin{pmatrix}
A_1 & & \\
& A_2 & \\
&&A_3
\end{pmatrix}, \qquad 
B_j = 
\begin{pmatrix}
b_{1j}I_N & & \\
& b_{2j}I_N & \\
& & b_{3j}I_N
\end{pmatrix}, \qquad j=1,2,3.
\end{equation*} 
Finally, as the FOM of the LVS, we obtain the following linear-quadratic system of ODEs 
\begin{equation}\label{odes}
\dot{\mathbf{u}} = \mathbf{F}(\mathbf{u}) := \underbrace{ A\mathbf{u} }_{linear} - \sum_{j=1}^3 B_j \underbrace{ (\mathbf{u}_{(j)}\circ\mathbf{u}) }_{quadratic}. 
\end{equation}

\subsection{Time discretization}

For time discretization, we employ Kahan's method, a nonstandard discretization technique \cite{Kahan93,Celledoni13} specifically developed for ODEs with quadratic nonlinearities.
For a general linear quadratic ODE $\dot{\mathbf{y}}=\mathbf{F}(\mathbf{y})=A\mathbf{y} + Q(\mathbf{y})$, the Kahan's method yields
\begin{equation}\label{genquad}
 \frac{{\mathbf y}^{n+1} - {\mathbf y}^n}{\Delta t} =   \frac{1}{2} A ({\mathbf y}^n + {\mathbf y}^{n+1}) + \widetilde{Q} ({\mathbf y}^n,{\mathbf y}^{n+1}),
\end{equation} 
with the bilinear form
$$
\widetilde{Q} ({\mathbf y}^n,{\mathbf y}^{n+1}) = \frac{1}{2} \left(Q({\mathbf y}^n+{\mathbf y}^{n+1}) - Q({\mathbf y}^n) -  Q ({\mathbf y}^{n+1}) \right).
$$
For the quadratic field $Q$, the scheme \eqref{genquad} is also equivalent to the equation \cite{Celledoni13}
$$
\left( I -\frac{\Delta t}{2} \mathbf{F}'(\mathbf{y}^n)\right)(\mathbf{y}^{n+1} - \mathbf{y}^n)  = \Delta t\mathbf{F}(\mathbf{y}^n),
$$
where $\mathbf{F}'$ denotes the Jacobian matrix of $\mathbf{F}$, and $\mathbf{y}^{n}=\mathbf{y}(t_{n})$ with the time instance $t_{n}=n\Delta t$, $n\geq 0$.
Kahan's method is a second-order, time-reversible, and linearly implicit scheme, meaning that each time step requires only a single Newton iteration.
In case of the FOM \eqref{odes}, the full discrete scheme reads as
\begin{align}\label{kahan}
\left( I_{3N} -\frac{\Delta t}{2} \mathbf{F}'(\mathbf{u}^n)\right)(\mathbf{u}^{n+1} - \mathbf{u}^n)  &= \Delta t\mathbf{F}(\mathbf{u}^n).
\end{align}

\section{Reduced order modeling} 
\label{sec:rom}

The construction of ROMs relies on the snapshot matrices obtained from the FOM \eqref{odes} of LVS, solved by the Kahan's scheme \eqref{kahan}.
The ROM solutions aim to approximate the FOM solutions within a low-dimensional linear subspace. 
For each component, we define the snapshot matrices $S_{j}$, which consist of the time evolution of the component $\mathbf{u}_j$, as follows
\begin{align*}
S_{j}&=\left[\mathbf{u}_j(t_1) \mathbf{u}_j(t_2) \cdots \mathbf{u}_j(t_K)\right]\in \mathbb{R}^{N\times K},& j=1,2,3.
\end{align*}

For either component $\mathbf{u}_j$, the low-dimensional linear subspace that approximately spans the column space of the related snapshot matrix $S_{j}$, and captures the important dynamics of the FOM, which is constructed as the space spanned by the column vectors of the basis matrices $V_{j,r}\in \mathbb{R}^{N\times r}$. Here, the number $r\ll N$ denotes the dimension of the reduced space, and is taken as the same for each component, but can also be chosen differently. Through the application of POD, the basis matrix $V_{j,r}$ is determined as the first $r$ left singular vectors corresponding to the dominant singular values from the SVD of the snapshot matrix $S_{j}$
\begin{align*}
 S_j &= V_{j,R_j}\Sigma_{j,R_j} W_{j,R_j}^T,& j=1,2,3,
\end{align*}
where $V_{j,R_j}\in {\mathbb R}^{N\times R_j}$ is the matrix of left singular vectors, $W_{j,R_j}\in {\mathbb R}^{K\times R_j}$ is the matrix of right singular vectors, $\Sigma_{j,R_j}\in {\mathbb R}^{R_j\times R_j}$ is the diagonal matrix whose diagonal entries are the singular values $\sigma_{j,1}\geq \cdots \geq \sigma_{j,R_j}\geq 0$, and $R_j$ is the rank of the snapshot matrix $S_{j}$.  
The POD basis minimizes the least squares errors
$$
\min_{V_{j,r}\in \mathbb{R}^{N\times r}}  \|S_j - V_{j,r} V_{j,r}^T S_j  \|_F^2 = \sum_{k=r+1}^{R_j} \sigma_{j,k}^2 ,
$$
where $\|\cdot\|_F$ is the Frobenius norm. The size of the POD basis is generally chosen by relative cumulative energy criteria which utilizes the singular values.

Since the low-dimensional linear subspace is spanned by the POD basis matrices $V_{j,r}$, we can write the low-rank solutions $\widehat{\mathbf{u}}_j(t)\in{\mathbb R}^{N}$ approximating the full-order solutions $\mathbf{u}_j(t)\in{\mathbb R}^{N}$ in the form
\begin{align}\label{lowappr}
\mathbf{u}_j(t)&\approx \widehat{\mathbf{u}}_j(t) = V_{j,r}\mathbf{w}_j(t), & j=1,2,3,
\end{align}
where $\mathbf{w}_j(t)\in{\mathbb R}^{r}$ is the vector of reduced coefficients.
Defining the $3r$-dimensional vectors $\mathbf{w}(t)$ and $\mathbf{w}_{(j)}(t)$ as
\begin{equation*}
\mathbf{w}(t) = 
\begin{pmatrix}
\mathbf{w}_1(t) \\
\mathbf{w}_2(t) \\
\mathbf{w}_3(t)
\end{pmatrix}, \qquad 
\mathbf{w}_{(j)}(t)  = 
\begin{pmatrix}
\mathbf{w}_j(t)  \\
\mathbf{w}_j(t)  \\
\mathbf{w}_j(t) 
\end{pmatrix}, \qquad j=1,2,3,
\end{equation*} 
and $3N\times 3r$-dimensional POD matrices $V_{r}$ and $V_{(j)}$ as
\begin{equation*}
V_r = 
\begin{pmatrix}
V_{1,r} & & \\
& V_{2,r} & \\
& & V_{3,r}
\end{pmatrix}, \qquad 
V_{(j)} = 
\begin{pmatrix}
V_{j,r} & & \\
& V_{j,r} & \\
& & V_{j,r}
\end{pmatrix}, \qquad j=1,2,3,
\end{equation*} 
we have that $\mathbf{u}(t)\approx \widehat{\mathbf{u}}(t) = V_r\mathbf{w}(t)$, where the reduced solution vector $\mathbf{w}(t)$ is computed from the following reduced system 
\begin{equation}\label{pod}
\dot{\mathbf{w}}(t) = A_r\mathbf{w}(t)  - \sum_{j=1}^3 V_r^TB_j ((V_{(j)}\mathbf{w}_{(j)}(t))\circ (V_r\mathbf{w}(t))). 
\end{equation}
The system \eqref{pod} is obtained by replacing $\mathbf{u}(t)\approx V_r\mathbf{w}(t)$ and projecting the FOM \eqref{odes} onto the low-dimensional subspace spanned by the POD basis $V_r$.
In the reduced system \eqref{pod}, $ A_r=V_r^TAV_r\in\mathbb{R}^{3r\times 3r}$ is the constant reduced matrix which can be precomputed. We call the reduced system of ODEs \eqref{pod} as the POD scheme, and like the FOM, it is also iterated in time by the Kahan's method.

\subsection{POD in tensor form}

Since the FOM \eqref{odes} contains nonlinear terms in the form of quadratic polynomials, it is desirable that its projection-based ROM maintains the same polynomial structure within the reduced subspace.
This can be handled using tensors, where each quadratic terms given by the element-wise product. The POD scheme \eqref{pod} is first converted into a form defined by tensor/Kronecker product.

Let $H\in\mathbb{R}^{3N\times (3N)^2}$ denote the matricized form of the tensor so that the identity $\mathbf{u}\circ\mathbf{v}=H(\mathbf{u}\otimes\mathbf{v})$ holds for any vectors $\mathbf{u},\mathbf{v}\in\mathbb{R}^{3N}$.
Then, by utilizing the properties of the Kronecker product, the quadratic term in \eqref{pod} can be reformulated as 
\begin{equation}\label{tensor}
\begin{aligned}
V_r^TB_j ((V_{(j)}\mathbf{w}_{(j)})\circ (V_r\mathbf{w})) &= V_r^TB_j H((V_{(j)}\mathbf{w}_{(j)})\otimes (V_r\mathbf{w})) \\
&= V_r^TB_j H(V_{(j)}\otimes V_r)(\mathbf{w}_{(j)}\otimes \mathbf{w})\\
&= H_{j,r}(\mathbf{w}_{(j)}\otimes \mathbf{w}),
\end{aligned}
\end{equation}
where the constant matrix $H_{j,r}=V_r^TB_j H(V_{(j)}\otimes V_r)\in\mathbb{R}^{3r\times (3r)^2}$ can be precomputable. 
By substituting the identity \eqref{tensor} into the system \eqref{pod}, we arrive-similarly to the FOM \eqref{odes}-at a linear-quadratic system of ODEs
\begin{equation}\label{tpod}
\dot{\mathbf{w}}(t) = A_r\mathbf{w}(t)  - \sum_{j=1}^3 H_{j,r}(\mathbf{w}_{(j)}(t)\otimes \mathbf{w}(t)). 
\end{equation}
We call the reduced system of ODEs \eqref{tpod} as the TPOD scheme. 
TPOD preserves the structure of FOM in the reduced space for PDEs and ODEs with quadratic nonlinearities. Additionally, because the online-offline computation is separated, the computational cost is lower than of the POD scheme \eqref{pod}.

The computational complexity in the offline stage is also decreased by calculating the precomputable matrix $H_{j,r}=V_r^TB_j H(V_{(j)}\otimes V_r)$, where the Kronecker product $V_{(j)}\otimes V_r$ has a computational complexity of ${\mathcal O}(r^2N^2)$. The matrix $H_{j,r}$ can be constructed without explicitly defining the matricized tensor $H$, which is represented in the MatLab notation as
\begin{align}\label{goyal}
	H_{j,r}&=V_r^TB_j H(V_{(j)}\otimes V_r)
	=V_r^TB_j
	\begin{pmatrix}
	V_{(j)}(1,:)\otimes V_r(1,:)\\
	\vdots\\
	V_{(j)}(3N,:)\otimes V_r(3N,:)
	\end{pmatrix},
	\end{align}
whose complexity is $ \mathcal{O}(Nr^3) $.
For any vectors $\textbf{u}$ and $\textbf{v}$, using the Kronecker product property
\begin{equation*}
(\text{vec}(\textbf{v\;u}^T)^T) = (u \otimes v)^T   = u^T \otimes v^T,
\end{equation*}	
 where $\text{vec}(\cdot)$ appends the columns of a matrix to form a vector, the matrix $ N_j:=H(V_{(j)}\otimes V_r)\in \mathbb{R}^{3N\times (3r)^2}$ in \eqref{goyal} can be computed as 
\begin{align*}
N_j(i,:)&=\left(\text{vec}\left(V_r(i,:)^T V_{(j)}(i,:)\right)\right)^T, & i=1,2,\ldots,3N.
\end{align*}
For an efficient calculation of the matrix $N_j=H(V_{(j)}\otimes V_r)$, we utilize \emph{multiprod} \cite{leva08mmm} which uses virtual array expansion to perform multiple matrix products. 
When the matrix $ V_{(j)}\in \mathbb{R}^{3N\times 3r} $ is reshaped into a three-dimensional tensor as $ V_{(j)} \in \mathbb{R}^{3N\times 1 \times 3r} $, the \emph{multiprod} operation is then performed between $V_r$ and $V_{(j)}$ along the second and third dimensions.
In this process, \emph{multiprod} implicitly treats $V_r$ as having a singleton third dimension, producing a three-dimensional tensor 
$$
\mathcal{N}_j:=multiprod(V_r,V_{(j)})\in\mathbb{R}^{3N\times 3r \times 3r}.
$$
As a result, equation \eqref{goyal} can be rewritten as $H_{j,r} = V_r^T B_j N_j^{(1)}$, where 
$N_j^{(1)}\in \mathbb{R}^{3N\times (3r)^2}$ denotes the matricized form of the tensor $\mathcal{N}_j$.

\section{Numerical results} 
\label{sec:num}

To show appearance of  segregation patterns, we consider the LV system \eqref{lv3} in a rectangular domain $\Omega =[0,2]\times[0,2]$  with 
the parameters $d_i=1.0e-2, \quad a_{ii}=1$  and $a_{ij}=3$ ($i\neq j$) \cite{Ei99}.
The initial conditions $u_i(x,0)$ are taken as small random perturbations around $0.143$.

\begin{figure*}[ht]
\centerline{\includegraphics[width=0.8\columnwidth]{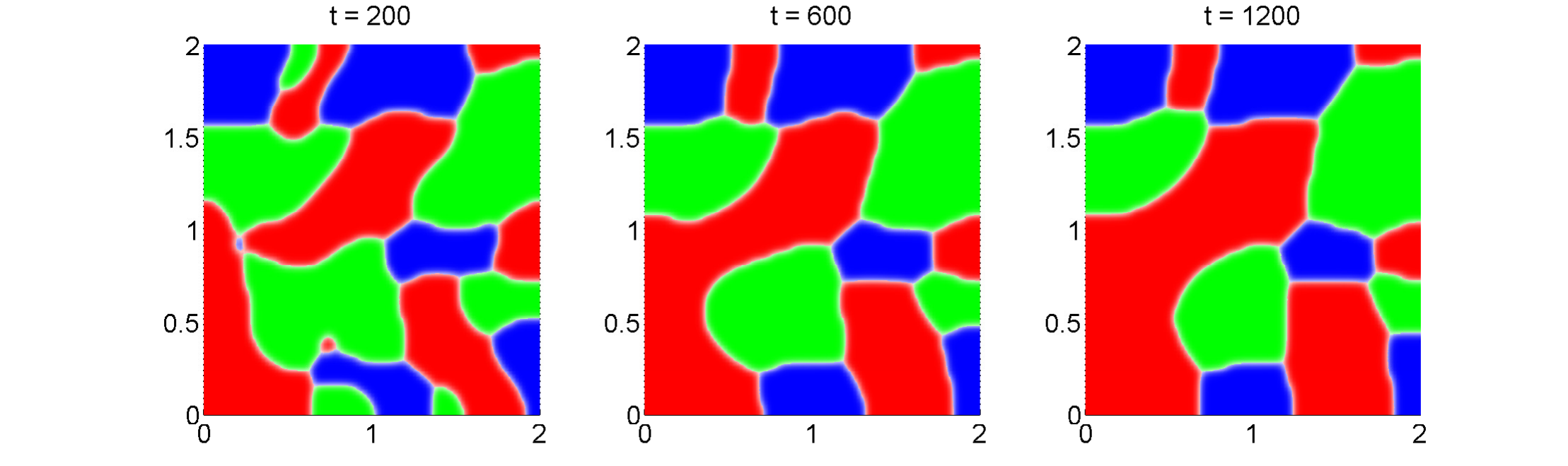}}
\caption{Phase separation at different times: $u$ (red), $v$ (green), $w$ blue.\label{solfom}}
\end{figure*}

In the initial stage, transition layers emerge, typically partitioning the domain $\Omega$ into sub-domains $\Omega_1$, $\Omega_2$, and $\Omega_3$. This suggests the onset of spatial segregation among three competing species, characterized by the formation of triple junctions.
We remark due to the random initial conditions, different patterns may emerge in long-time computation for each run. In Figure~\ref{solfom}, the appearance of some of the solution profiles at times $t=200, 600, 1200$ are presented.

\begin{figure*}[ht]
\centerline{\includegraphics[width=0.5\columnwidth]{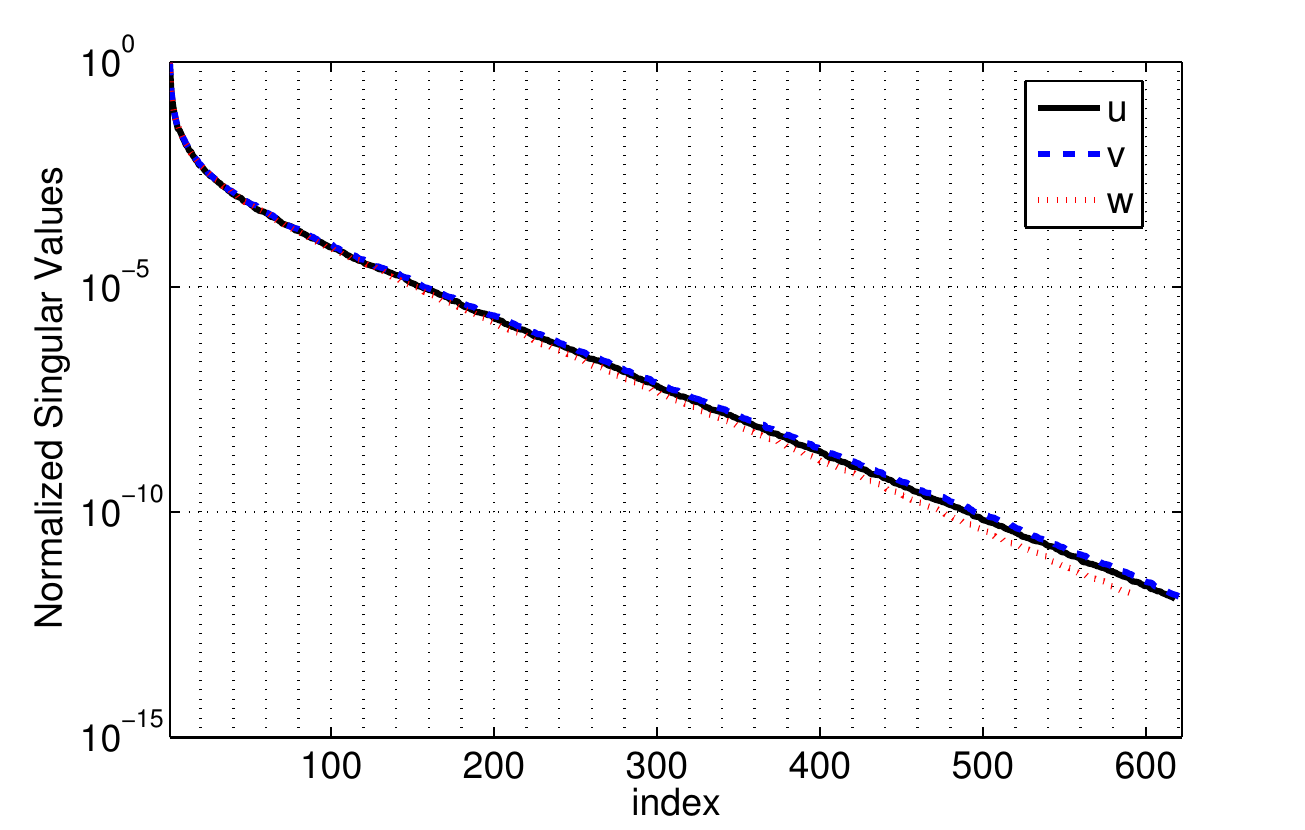}}
\caption{Decay of singular values.\label{sing}}
\end{figure*}

The singular values decay slowly in Figure~\ref{sing}, as for other pattern forming systems  \cite{Alla23,Bozzini21}, FitzHugh-Nagumo equation \cite{Karasozen17}, and Shigesada-Kawasaki-Teramoto cross-diffusion system \cite{Karasozen21}. 
This leads to requirement of large number of POD modes to represent the FOM accurately, but effects the computational efficiency of the ROMs. 
Therefore, we have restricted the number of POD modes to $r=20$ which is sufficient to resolve the patterns accurately by the ROMs in long-time computation.

\begin{figure*}[ht]
\centerline{\includegraphics[width=0.9\columnwidth]{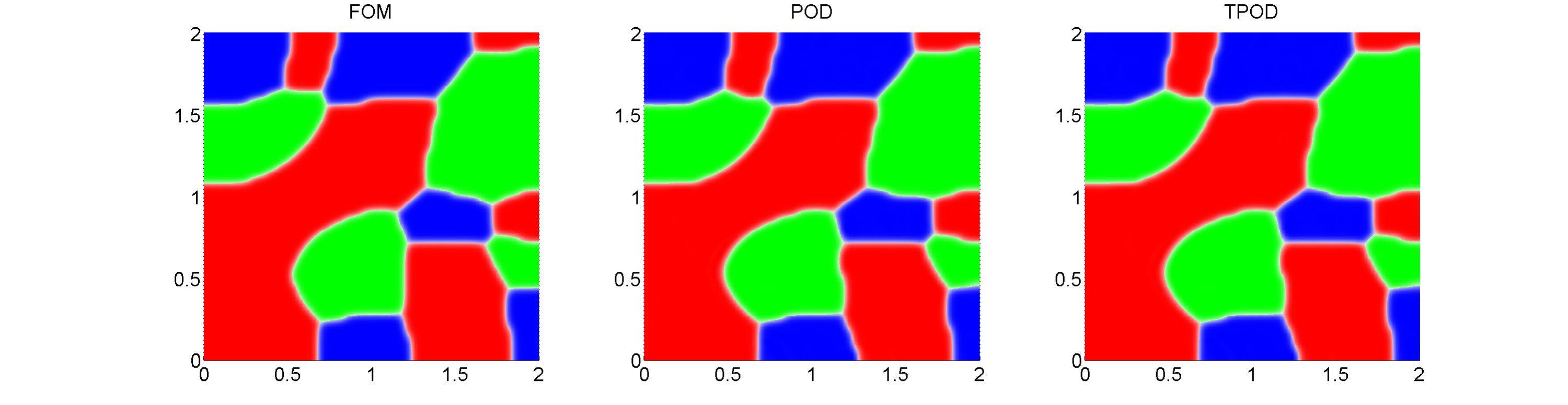}}
\caption{Phase separation by FOM/ROM solutions at final time: $u$ (red), $v$ (green), $w$ blue.\label{pat}}
\end{figure*}

As the first course, we investigate the performance of the ROMs using $r=20$ POD modes.
The segregation patterns of the ROMs are very close to those of the FOM in  Figure~\ref{pat}, which indicates the sufficiently well-behavior of the POD and TPOD schemes.

\begin{figure*}[ht]
\centerline{\includegraphics[width=0.5\columnwidth]{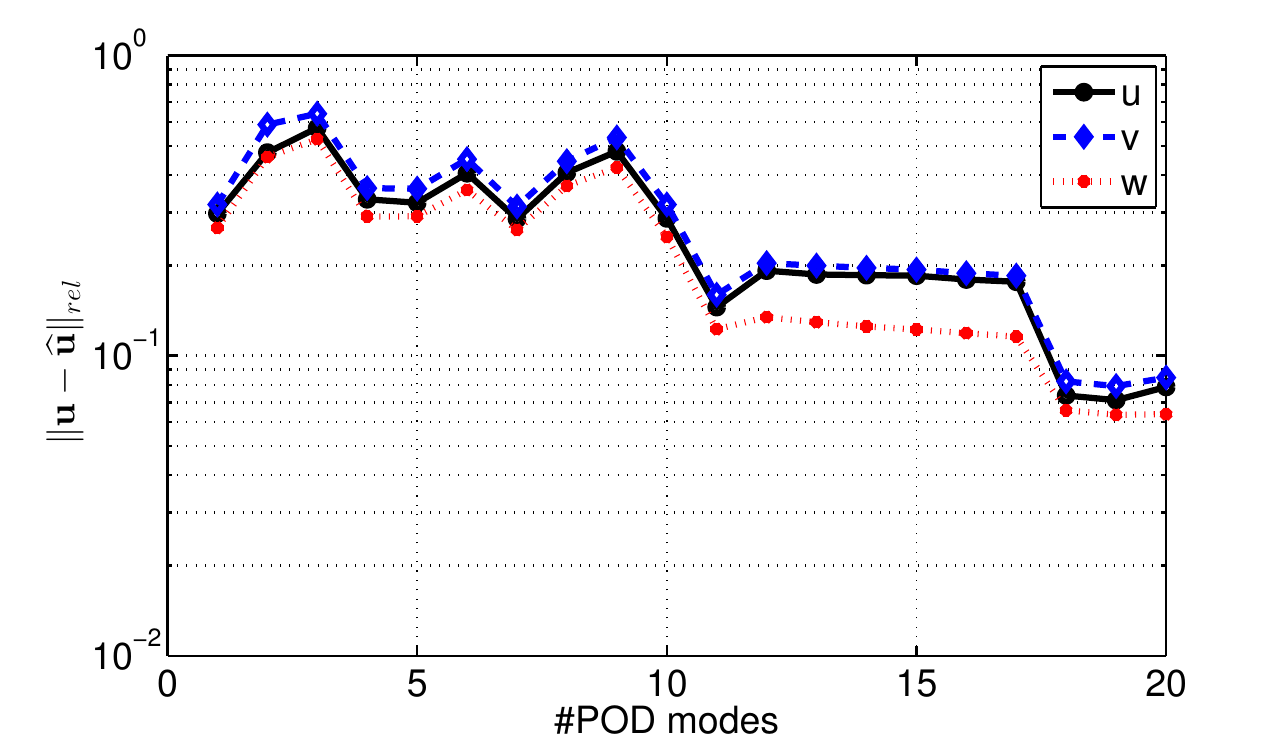}}
\caption{Relative FOM-ROM errors.\label{L2L2}}
\end{figure*}

For the accuracy of the ROMs, the time-averaged relative $L^2$-errors between the FOM and ROM solutions are computed for each state variable ${\mathbf u}= {\mathbf u_1}, {\mathbf u_2},{\mathbf u_3} $
\begin{align}\label{relerr}
\|\mathbf{u}-\widehat{\mathbf   u}\|_{rel}=\frac{1}{K}\sum_{k=1}^{K}\|\mathbf{u}-\widehat{\mathbf   u}\|_{rel,k},
\end{align}
where $\widehat{\mathbf  u}$ being the reduced approximation to ${\mathbf u}$, $u^k_{i,j}\approx u(x_i,y_j,t_k)$, and 
\begin{align}\label{relerrk}
\|\mathbf{u}-\widehat{\mathbf   u}\|_{rel,k}=\frac{\|{\mathbf  u}^k-\widehat{\mathbf  u}^k\|_{L^2}}{\|{\mathbf  u}^k\|_{L^2}}, \quad  \|{\mathbf  u}^k\|_{L^2}^2=\sum_{i=1}^n\sum_{j=1}^n (u^k_{i,j})^2\Delta x\Delta y.
\end{align}
The relative FOM-ROM errors in Figure~\ref{L2L2}, tends to decrease with increasing POD modes.

For the computational efficiency gained by TPOD, we present in Table~\ref{cpu} the wall clock times in seconds for FOM and ROMs computations.  
It is important to note that the wall clock time for ROMs computations includes not only the online simulation time but also the calculation of the POD basis functions for POD, and both the calculation of the POD basis functions and the computation of the tensorized matrices for TPOD. The speed-up factors presented in Table~\ref{cpu} are determined by comparing the wall clock time needed to obtain the FOM solutions with the one for the ROM solutions. These factors highlight the computational efficiency of the ROM by exploiting tensor techniques.

\begin{table}
\caption{Wall-clock time and speed up factors with the reduced dimension $r=20$ \label{cpu}}
\centering
\begin{tabular}{l|l|l|l}
\hline
 &  \textbf{FOM} & \textbf{POD}  & \textbf{TPOD} \\
\hline
Total Computation time & 565.1 & 358.4  & 50.1  \\
\hline
 Speed up              &    -    & 1.58 & 11.28 \\
\hline
\end{tabular}
\end{table}

\begin{figure*}[ht]
\centerline{\includegraphics[width=0.5\columnwidth]{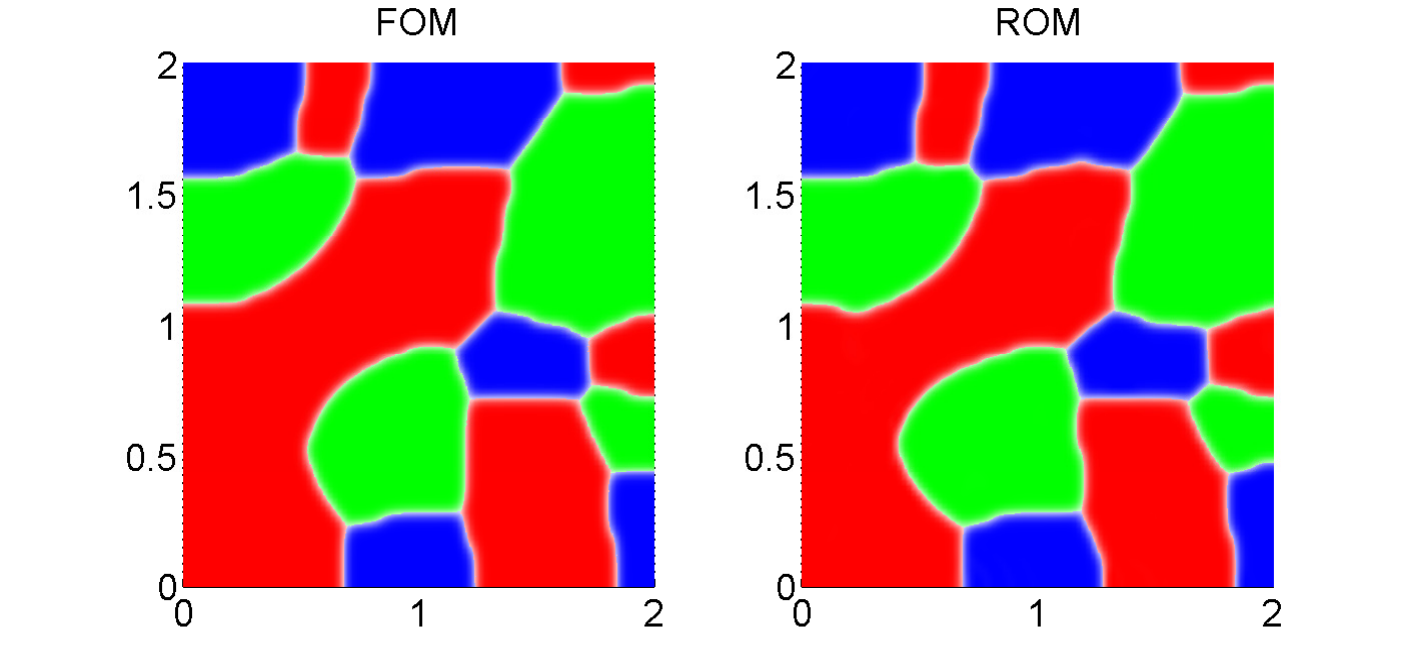}\includegraphics[width=0.36\columnwidth]{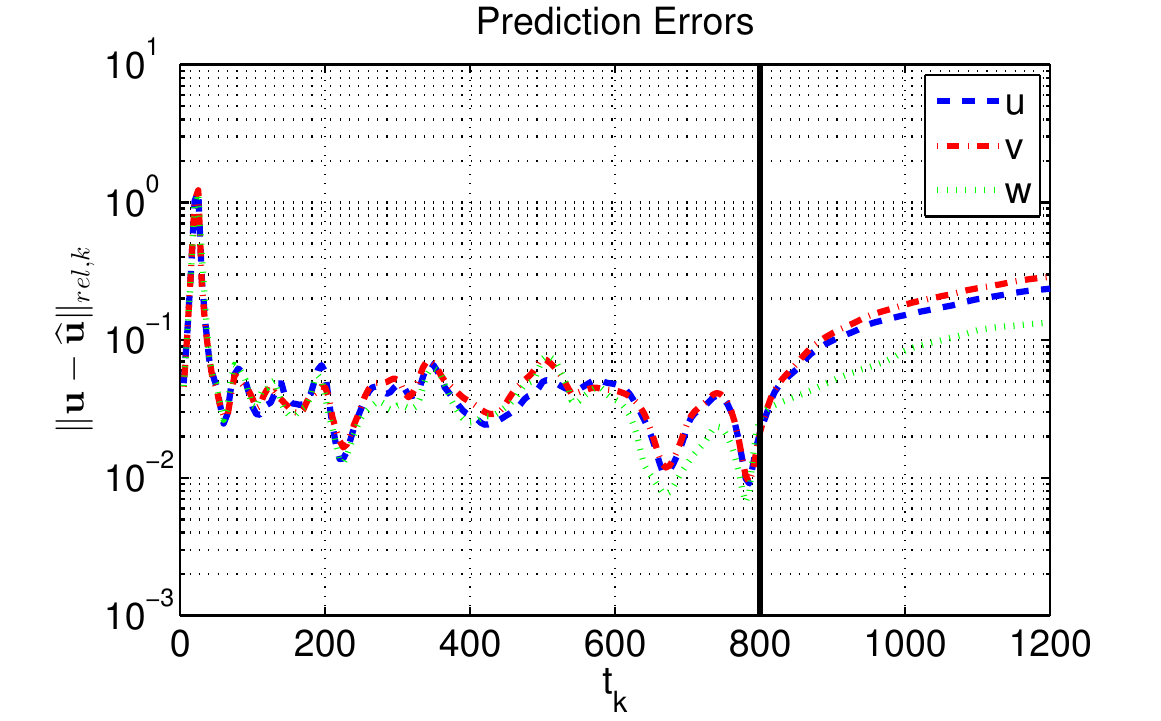}}
\centerline{\includegraphics[width=0.5\columnwidth]{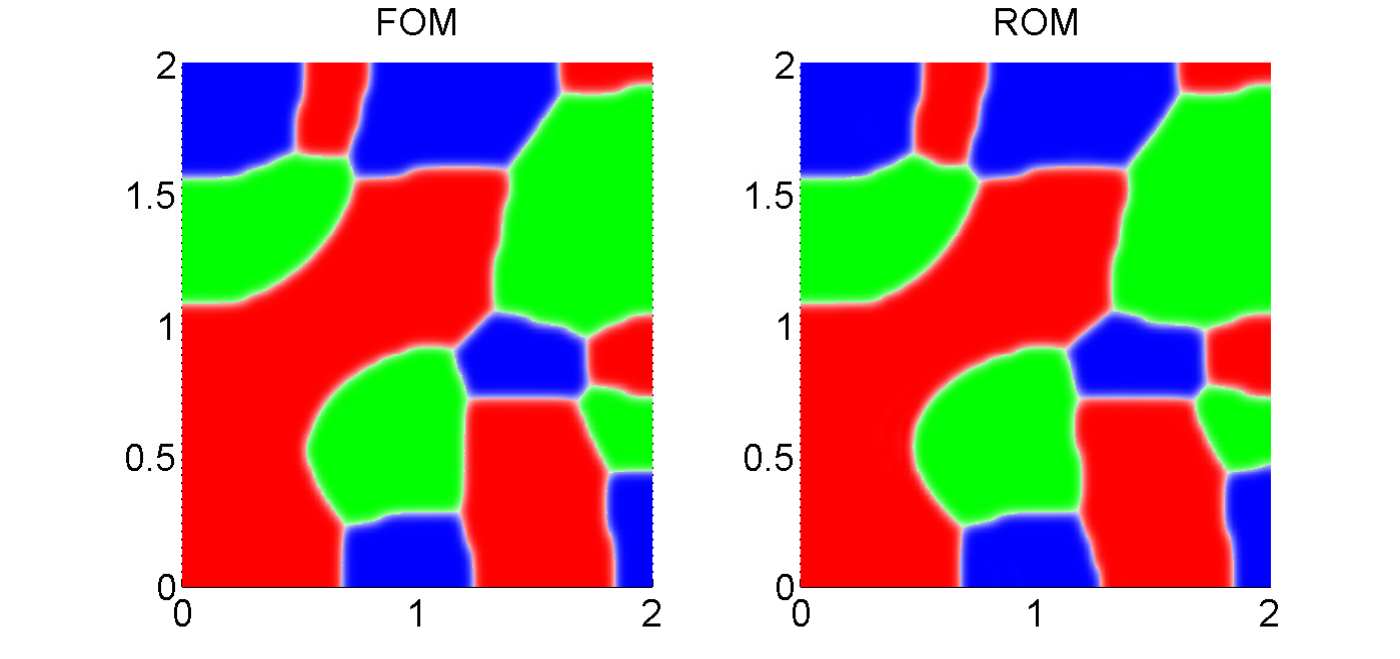}\includegraphics[width=0.36\columnwidth]{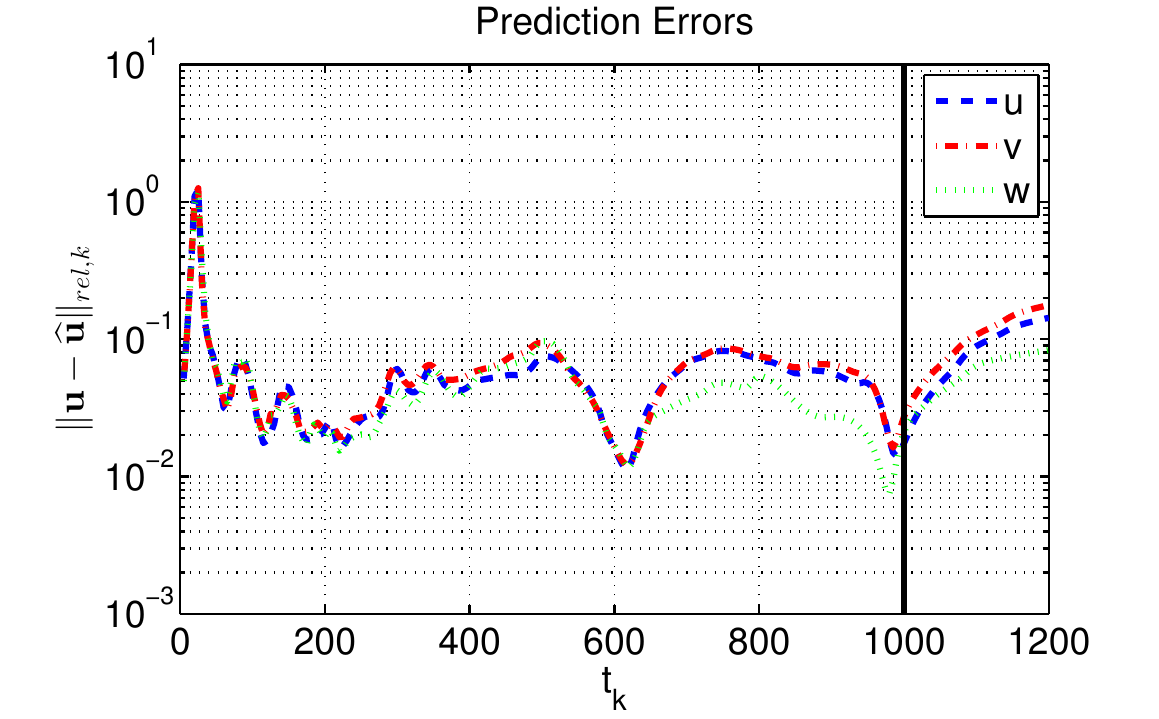}}
\caption{Predicted phase separation at final time and relative $L^2$-errors: the vertical line indicates the time which separates the learning and prediction phases.\label{pred}}
\end{figure*}

As the last course, we demand to utilize the ROM (TPOD) to predict the solution patterns after a given time instance $t_L<T$. 
More precisely, we are about to find the POD basis modes learned by the FOM solution snapshots on the interval $[0,t_L]$ (learning phase), and then using the ROM scheme based on these modes, we expect to predict the patterns for $t\in [t_L,T]$ (prediction phase). 
In Figure~\ref{pred}, the phase separation of FOM/ROM patterns at final time $T=1200$, and the relative $L^2$-errors \eqref{relerrk} between FOM and ROM solutions along the time trajectory are shown.
In Figure~\ref{pred}, top, the ROM solutions are obtained with the learning time $t_L=800$, and for  $t_L=1000$ at the bottom.
The truncation errors decrease with the learning time $t_L$ when $t\in [800,100]$, resulting in a more accurate prediction, because the ROMs can better capture the dynamics by learning from larger data sets. However, we see that the errors at the final time is about $10^{-1}$ in either case.

\section{Conclusions}
\label{sec:con}

In this work, we developed and analyzed ROM techniques for the efficient approximation of spatio-temporal patterns in three-species diffusive Lotka-Volterra systems. The FOM, constructed via the SEM and time-integrated using Kahan's linearly implicit scheme, successfully captures the emergence of complex spatial segregation patterns typical of competitive population dynamics. To address the computational demands of long-time simulations on fine spatial and temporal grids, we applied POD to generate low-dimensional ROMs. Furthermore, we introduced TPOD formulation that exploits the quadratic structure of the underlying dynamics and enables a clear offline-online computational separation. Our numerical results demonstrate that TPOD significantly reduces computational time-achieving a speedup factor more than an order-of-magnitude, while maintaining accuracy in reproducing the full model dynamics.

We also explored the predictive capabilities of the ROMs by training them on partial time intervals and extending the simulation beyond the training horizon. The ROMs were shown to accurately capture the future evolution of the patterns, provided sufficient training data.
Overall, the combination of linearly implicit time integration, spectral discretization, and tensor-based reduced order modeling proves to be an effective and scalable strategy for simulating complex reaction-diffusion systems. Future work may consider extending these techniques to more general nonlinear systems, adaptive basis enrichment, and integration with data-driven modeling frameworks.


\end{document}